\documentstyle[12pt,psfig]{article}
\input xy
\xyoption{all}
\begin{document}

\hfill DUKE-CGTP-2001-14

\hfill math.DG/0110156

\vspace*{1.5in}

\begin{center}

{\Large\bf Discrete Torsion, Quotient Stacks,}

{\Large\bf and String Orbifolds}

\vspace{0.5in}

{\bf Lecture at Wisconsin-Madison}

{\bf Workshop on Mathematical Aspects of Orbifold String Theory}

{\bf May 4-8, 2001}

\vspace{0.5in}

Eric Sharpe \\
Department of Physics \\
Box 90305 \\
Duke University \\
Durham, NC  27708 \\
{\tt ersharpe@cgtp.duke.edu} \\

\end{center}

%

\newpage

\tableofcontents

\newpage

\section{Introduction}

In this talk we shall describe two separate topics:  discrete torsion
(a previously-mysterious degree of freedom in string orbifolds,
of which we shall give a reasonably complete understanding),
and the relation between quotient stacks and string orbifolds
(where, unlike discrete torsion, we shall only set up basics,
and will emphasize that much work remains to be done).

Discrete torsion is a degree of freedom that appears in describing
string orbifolds. 
Historically discrete torsion has been considered very mysterious.
However, in the first part of this note we shall outline recent work 
\cite{dt3,cdt,hdt,dtrev} that de-mystifies it.
To be brief, we shall argue that discrete torsion is the choice
of equivariant structure ({\it i.e.}, orbifold group action)
on the $B$ field, and from this derive the classification by
$H^2(G, U(1))$, Vafa's twisted sector phases, Douglas's projectivized
D-brane actions ({\it i.e.}, projectivized equivariant K-theory), and
analogues for other tensor field potentials.  In a nutshell,
discrete torsion has a natural understanding that has nothing to
do with conformal field theory, Riemann surfaces, or any other
baggage of perturbative string theory, and we shall outline
this understanding.

In the second part of the talk we shall discuss the relation between
string orbifolds and quotient stacks, and in particular we shall
outline how a string orbifold is precisely a sigma model on a quotient
stack, a description that clarifies the physics
of string orbifolds.  For physicists, this notion is a 
radical conceptual shift:
although string orbifolds are described in terms of group actions on
covers, physicists have historically assumed that this was merely
scaffolding.
Physicists speak of string orbifolds as describing strings on
quotient spaces decorated with some sort of quantum or `stringy'
behavior at the singularities; for example, 
physicists often speak of string orbifolds as describing
strings on quotient spaces suitably 
decorated with $B$ fields \cite{edstrings95,paulz2,kw},
or as strings on some resolution of the quotient space
(because of massless moduli which can often be interpreted as
K\"ahler moduli).
In particular, no physicist has 
ever claimed\footnote{After all, one can use group actions on covers
to describe quotient spaces as well as quotient stacks, and quotient spaces
are far less exotic.  
} that string orbifolds describe strings on any
sort of stack, or (equivalently)
that any formal geometric meaning assigned to
the group-actions-on-covers scaffolding had any physical relevance.
To make matters even more confusing, in practice physicists often 
implicitly assume that string orbifolds describe strings
on quotient spaces, and ignore any quantum effects
at singularities.  For example, string moduli spaces are constructed
around this assumption, and they play an important role in
understanding many of the duality symmetries that have been
of interest to physicists over the last decade.


For mathematicians who are acquainted with quotient stacks,
the idea that string orbifold conformal field theories (CFT's) coincide with
CFT's for strings 
compactified on quotient stacks
surely seems much more natural.  After all, among other things, quotient stacks
are an overcomplicated way to describe group actions on covers,
the language used in string orbifolds.  This formal similarity
might even lead someone who was not acquainted with the physics literature
to claim that they `know' string orbifolds are strings on 
quotient stacks.  However, for such a statement to be true
implies that string orbifold CFT's coincide with 
CFT's for strings on quotient stacks, and even assuming that the
notion of compactification on stacks is sensible, 
there is a tremendous amount of work that must be done to justify this.
Put another way, such a statement implicitly assumes that the
extra structure of a `generalized space,' as possessed by a stack,
is physically relevant.
Any competent physicist would observe that not only
is the physics lore
apparently contradictory\footnote{After all,
if string orbifolds really do describe strings on resolutions,
then they cannot possibly also describe strings on quotient stacks --
clearly, a choice must be made.}, but one would need to understand
string compactification on stacks before such a statement could
really be justified, and furthermore, such a statement fails several
basic physical consistency conditions.
In more detail:
\begin{itemize}
\item Before one can claim that string orbifold CFT's coincide with
CFT's for strings 
compactified on
{\it quotient} stacks, one must first check whether the notion of string 
compactification on stacks is even sensible,
something that was not considered by physicists at all
until very recently \cite{qstx}.
One way to do this (which we shall describe the first stages of)
is to first write down the classical action for a sigma model on
a stack, and understand basic consequences of that notion,
such as the massless spectrum of that sigma model.
To be certain that such a classical action can be consistently
quantized, there are global considerations that must be taken
into account.   
One must resolve apparent physical contradictions, such as the fact
that the massless spectrum of such sigma models is {\it not} given
by cohomology of the target.
Countless questions, ranging from ``how does one make sense of anomalies
in this context'' (and other nontrivial global issues)
to ``can one do QFT on a stack, now viewed
as spacetime,'' must be answered
before the matter can be considered to be completely settled.
After one has settled at least some of these issues,
one can then check explicitly whether, in fact, CFT's for strings
compactified on quotient stacks really do coincide with CFT's for
string orbifolds.
\item The statement that a string orbifold CFT coincides with the CFT
for a string compactified on a quotient stack (assuming this is
a sensible notion) has nontrivial physical implications,
which must be checked for consistency.
Put another way, one cannot consistently consider
string orbifolds in isolation from the rest of string theory.
For example, in constructing moduli spaces of
string vacua (used to justify items from mirror symmetry to
string/string duality)
physicists have always assumed that the deformation theory
of a string orbifold is the same as that of a quotient space.
If string orbifolds do not describe strings on quotient spaces,
then one must explain how the deformation theory arguments used
by physicists can possibly have been consistent.
\end{itemize}
In short, not only does the idea that string orbifolds coincide
with string compactification on quotient stacks naively contradict
the physics lore, but a tremendous amount of very basic work must be done
to begin to justify such claims, and such claims even appear
to yield physical contradictions.  We shall describe some of the basic
work needed to justify such claims, and resolve some (but not all)
of the contradictions, but much work remains to be done.
In particular, it must be said that, at present, there is still a good chance
that string orbifolds do {\it not} describe strings on quotient stacks.

Another question one might ask is, why bother?
If one is only using quotient stacks as a highly overcomplicated 
means of describing group actions on covers, then there is hardly a point.
However, we shall argue later that the idea that a string orbifold CFT
coincides with the CFT for a string compactified on a quotient stack
(if indeed this is a sensible notion) has highly nontrivial physical
implications.  One aspect is that this gives a new geometric way
of understanding certain 
physical properties of string orbifolds.
Another aspect, as mentioned above,
is that this calls into question the arguments physicists have used to
construct moduli spaces of string vacua, essential to understand
string dualities of all types.


As we have previously written about quotient stacks for
a physics audience \cite{qstx}, here we shall 
speak to a mathematics audience.
We shall describe some of the basics needed by physicists to make sense of
the notion that string orbifolds describe strings compactifieid
on quotient stacks,
and shall outline the work that remains to be done before
quotient stacks can be universally accepted in the physics
community as being genuinely relevant to string orbifolds.


\section{Lightning review of string orbifolds}   \label{lightning}

Before describing either discrete torsion or the relationship
between string orbifolds and quotient stacks,
we shall take a few moments to review string orbifolds.
A string orbifold \cite{dhvw1,dhvw2} is simply a sigma model with the action of
a discrete group on the target space gauged.
Since a string orbifold is a gauged string sigma model,
let us take a moment to review sigma models.

\begin{figure}
\centerline{\psfig{file=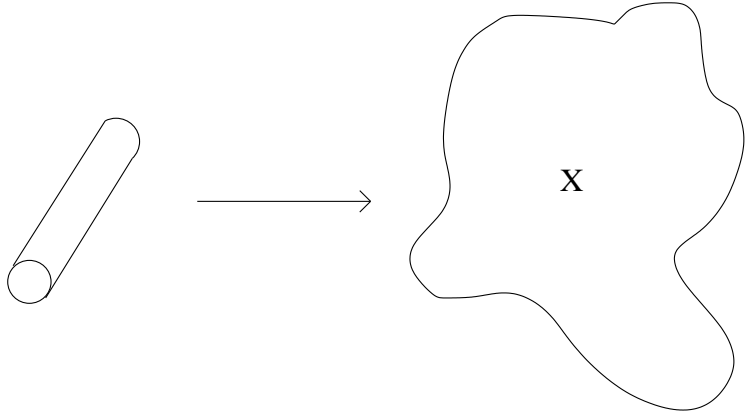,width=4in}}
\caption{\label{gensig} String sigma models sum over ``worldsheets''
swept out by strings in $X$, as shown for a free string,
whose worldsheet is the cylinder $S^1 \times {\bf R}$.}
\end{figure}

A sigma model with target space $X$ and base space\footnote{
For standard sigma models, $Y$ is assumed a manifold with either
a Lorentzian or Riemannian metric.  Depending upon whether one
wants to describe all of spacetime, or just a factor in a spacetime
of the form ${\bf R}^4 \times Y$, for example, one can consider
either case of Lorentzian or Riemannian signature.
Analogous issues arise for worldsheet metrics.  We shall ignore
this issue in the remainder of this section.} $Y$ 
is a weighted sum over maps $Y \longrightarrow X$,
as schematically illustrated for two-dimensional $Y$ 
in figure~\ref{gensig}.
(The `sum' in question is a sum in the sense of path integrals.)

For example, a sigma model describing the propagation of a point particle
on $X$ is a sum over maps from a one-dimensional $Y$
(one-dimensional because it encodes the `worldline' of a point particle
moving in spacetime) into $X$, weighted by $\exp( i S )$,
where $S$ is known as the classical action and has the form
\begin{equation}
S \: \sim \: \int dt (\phi^* G_{\mu \nu}) \frac{ d \phi^{\mu} }{ dt }
\frac{ d \phi^{\nu} }{ dt } \: + \: \cdots
\end{equation}
Such a sigma model is one description of the quantum mechanics of
a point particle on $X$ \cite{feynman}.

For another example, a sigma model describing the propagation of
a string on $X$, as illustrated in figure~\ref{gensig},
is a sum over maps from a two-dimensional $Y$
(two-dimensional because it encodes the path swept out by
the string over time), known as the worldsheet, into $X$,
weighted by $\exp(i S)$, where $S$ is known as the classical action
and has the form
\begin{equation}   \label{stdclassact}
S \: \sim \: \int d^2 \sigma \left( \phi^* G_{\mu \nu} \right)
h^{\alpha \beta}
\frac{ \partial \phi^{\mu} }{ \partial \sigma^{\alpha} } 
\frac{ \partial \phi^{\nu} }{ \partial \sigma^{\beta} } 
\: + \: \cdots
\end{equation}
Just as a sigma model for a point particle describes the quantum
mechanics of a point particle on $X$, a sigma model for a string
describes the `stringy quantum mechanics' of a string on $X$.
Formally one can continue to write sigma models in higher dimensions,
but above two dimensions they become less well-behaved.

Now, a physical orbifold is obtained by starting with a sigma
model on some space $X$, and `gauging' the action of a discrete group
$G$ on $X$.  To `gauge' the action of a discrete group means that
fields in the sigma model differing by the action of $G$ should be
identified -- we impose an equivalence relation on the `field space'
that our `path integral' integrates over, and integrate over equivalence
classes.  (It should be noted that the gauging we describe is a process 
performed
by physicists to build one class of physical theories from another class;
gauging a symmetry is not a physical process, but rather a process
performed by physicists.)

\begin{figure}
\centerline{\psfig{file=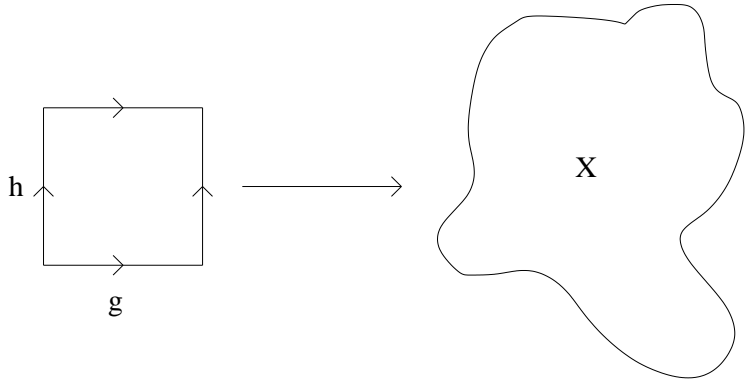,width=4in}}
\caption{\label{figsig1} 
A contribution to the $(g,h)$ twisted sector of a string orbifold $[X/G]$ on
$T^2$}
\end{figure}

In practice, what effect does this `gauging' have on a sigma model?
The answer is that in the new physical theory obtained
by gauging, fields on the base space $Y$, such as the
map into $X$, need no longer be well-defined over $Y$, but only
well-defined up to the action of $X$ -- we are allowed to introduce
branch cuts.  Also, as part of the gauging, we are required to sum
over all possible choices of branch cuts.  For example, if $Y = T^2$,
then instead of summing over maps $T^2 \longrightarrow X$,
we sum over maps as illustrated in figure~\ref{figsig1}, where branch
cuts have been introduced on the $T^2$.  Maps with nontrivial
branch cuts are known as contributions to a `twisted sector.'

More formally, each twisted sector contribution should be thought
of as, a $G$-equivariant map from the total
space of a principal $G$-bundle on $Y$ to $X$, restricted to a
particular lift of $Y$ to the total space of the bundle.
Such a lift introduces branch cuts, as shown in figure~\ref{figsig1}.
Then, our path integral sum is a sum over equivalence classes of
bundles and $G$-equivariant maps.

If we let $Z_{(g,h)}$ denote a (path-integral-type) sum over contributions
to the $(g,h)$ twisted sector on $T^2$, weighted as above by\footnote{Note
that we can only consider group actions for which the action $S$ is
well-defined on the twisted sector illustrated in figure~\ref{figsig1}.
This one of several constraints on possible choices of groups and
group actions.}
$\exp(i S)$, then naively we are led to believe that the path integral
sum (known as a `partition function') 
for a string orbifold on $T^2$ has the form
\begin{displaymath}
Z(T^2) \: = \: \sum_{\begin{array}{c} \scriptstyle{ g,h \in G } \\ 
\scriptstyle{ gh=hg}  \end{array}
 } Z_{(g,h)} 
\end{displaymath}

The expression above is almost correct, except for one
small subtlety.  If we sum over {\it all} possible twisted-sector
maps of the form
illustrated in figure~\ref{figsig1}, then we actually overcount
by $|G|$.  After all, as mentioned earlier, we only wish to sum
over equivalence classes of bundles and $G$-equivariant maps from
the total space of the bundle to $X$.  For each such bundle and
$G$-equivariant map, there are $|G|$ twisted-sector maps,
corresponding to $|G|$ distinct lifts, so summing over twisted sector
maps overcounts by $|G|$, and we find that the correct expression
for the path integral sum ({\it i.e.} the partition function)
is given by
\begin{equation}    \label{zt2}
Z(T^2) \: = \: \frac{1}{|G|}
\sum_{\begin{array}{c} \scriptstyle{ g,h \in G } \\
\scriptstyle{ gh=hg}  \end{array}
 } Z_{(g,h)}
\end{equation}

Now, although we have been talking about string orbifolds,
and string orbifolds on $T^2$, the same remarks apply not only
to other Riemann surfaces, but to $Y$'s of any dimension.
Again, for $Y$ of dimension greater than two, the notion of
a sigma model is not well-defined beyond the classical limit,
but one can write down such gauged sigma models classically in any
number of dimensions.

A few notes on expression~(\ref{zt2}) are in order.
\begin{enumerate}
\item First, note that $Z_{(1,1)}$ is the same as the partition function
for the original (ungauged) sigma model into $X$.
\item Second, note that once one introduces some nontrivial twisted sectors
into the theory on $T^2$,
consistency with modular invariance forces one to sum over
all possible twisted sectors.  After all, under the
transformation $\tau \mapsto \tau + 1$, a twisted sector
$(g,h) \mapsto (gh,h)$, so clearly modular transformations mix
twisted sectors.   
\item Third, some remarks on the Hamiltonian description of
orbifolds are in order.  From expression~(\ref{zt2}), as a string
propagates around the loop, it comes back to itself, but meets
the operator
\begin{displaymath}
\frac{1}{|G|} \sum_g g
\end{displaymath}
This is a projection operator, and it projects onto $G$-invariant
states -- only $G$-invariant string states are allowed to propagate.
\end{enumerate}

What can we do with this physical theory?
Given a string sigma model on a target space $X$,
the Euler characteristic of the target space $X$ can be obtained by
evaluating the partition function for the theory on $T^2$
in a limit of the worldsheet metric.  
Now, our `string orbifold' theory is no longer a sigma model
on any particular space, however one can play the same formal game
to recover the so-called `stringy orbifold Euler characteristic'
\cite{dhvw1,dhvw2,hirzhofer,atsegal}
\begin{eqnarray}
\chi_G(X) & = & \frac{1}{|G|}  
\sum_{\begin{array}{c} \scriptstyle{ g,h \in G } \\ 
\scriptstyle{ gh=hg}  \end{array}  }
e\left(X^{<g,h>}\right) \\
& = & \sum_{[g]} e\left( X^g / C(g) \right)   \label{hhoec}
\end{eqnarray}
where $X^{<g,h>}$ denotes the subset of $X$ invariant under both
$g, h \in G$, $e(X)$ denotes the Euler characteristic of $X$,
$[g]$ denotes the conjugacy class of $g \in G$, and
$C(g)$ denotes the centralizer of $g \in G$.

In the special cases when the quotient space $X/G$ admits
a crepant resolution, the orbifold Euler characteristic above
agrees\footnote{Moreover, the Euler characteristic is independent
of the choice of crepant resolution, when more than one exists.} 
with the Euler characteristic of the crepant resolution.
More generally, the physical theory has massless modes which
correspond to K\"ahler moduli (when the quotient space admits
crepant resolutions)i, and the conformal field theory behaves
as though it describes a string on a smooth space.  
For such reasons, physicists have 
historically often claimed that string orbifolds appeared to be
describing strings on some sort of resolution of the quotient
space.

More generally, how do physicists interpret string orbifolds?
String orbifolds were originally created in an attempt to
describe strings on quotient spaces.  However, they are not
sigma models on quotient spaces.
Also, as described above, they have certain physical properties which suggest
that they might have some sort of interpretation as strings
on resolutions (though in general resolutions need neither exist
nor be unique).  More recently, it has been suggested
\cite{edstrings95,paulz2,kw} that there is an accompanying 
$B$ field decoration,
which ties into other physical characteristics of string orbifold
conformal field theories.
To make matters even more confusing, in practice physicists often
implicitly assume that string orbifolds describe strings on quotient
spaces\footnote{Also see, for example, \cite{eva1,ach1} for a random
sampling of some recent prominent physics papers making this assumption in
different contexts.}.  
For example, string moduli spaces are constructed with the
implicit assumption that the deformation theory of a string orbifold
coincides with that of a quotient space.

One interpretation that most physicists would agree upon is that
string orbifolds describe strings on quotient spaces, but with
some sort of `stringy' behavior located at the singularities,
which has the effect of somehow resolving the singularities.
Since, in a sigma model, massless modes are identifiable with
cohomology of the target space, physicists have believed
that twist fields should have some understanding
as some unknown cohomology of the quotient space\footnote{To be
contrasted with a cohomology constructed from group actions on covers.
Although string orbifolds are phrased in terms of group actions on covers,
the point here is that they naively seem to predict the existence
of a cohomology theory directly on quotient spaces.},
referred to as
`orbifold cohomology.'  
It was hypothesized that knowledge of such an orbifold cohomology,
of the form suggested by physics,
would shed light on the physics underlying string orbifold
CFT's, by giving a better understanding of the stringy phenomena
taking place at the quotient singularities.

\section{Discrete torsion}

In the beginning we mentioned that discrete torsion is a degree
of freedom associated with string orbifolds.  How does it enter?

Discrete torsion was originally discovered in the following
fashion \cite{vafa1}.  Start with a string orbifold partition function
on, say, $T^2$, as we discussed earlier (equation~(\ref{zt2})):
\begin{displaymath}
Z(T^2) \: = \: \frac{1}{|G|} 
\sum_{\begin{array}{c} \scriptstyle{ g,h \in G } \\ 
\scriptstyle{ gh=hg}  \end{array}
 } Z_{(g,h)} 
\end{displaymath}
Now, weight the twisted sectors by phases, to obtain a new
partition function:
\begin{displaymath}
Z'(T^2) \: = \:  \frac{1}{|G|}
\sum_{\begin{array}{c} \scriptstyle{ g,h \in G } \\ 
\scriptstyle{ gh=hg} \end{array} }
\epsilon(g,h) \, Z_{(g,h)}
\end{displaymath}
$Z'$ is now the partition function of a theory ``with discrete torsion.''
The phases $\epsilon(g,h)$ are heavily constrained by internal consistency
conditions.  After one does some work, one finds that one solution
of the constraints is given by
\begin{equation}    \label{vafa1loop}
\epsilon(g,h) \: = \: \frac{  \omega(g,h)  }{   \omega(h,g)  }
\end{equation}
where the $\omega(g,h)$ are 2-cocycle (inhomogeneous) representatives 
of a class
in the group cohomology group\footnote{This group cohomology group
is defined with trivial action on the coefficients.  The same will
be true of all group cohomology referenced in this lecture.}
$H^2(G, U(1))$.
(For those readers who do not have group cohomology at their
fingertips, this just means that the $\omega$ are maps
$G \times G \rightarrow U(1)$, obeying the cocycle condition 
\begin{equation}    \label{2cocycle}
\omega(g_1 g_2, g_3) \, \omega(g_1, g_2) \: = \:
\omega(g_1, g_2 g_3) \, \omega(g_2, g_3)
\end{equation}
and with coboundaries defined by
\begin{displaymath}
\omega(g,h) \: \sim \: \omega'(g,h) \: \equiv \: f(gh) \,
\omega(g,h) \, f(g)^{-1} \, f(h)^{-1}
\end{displaymath}
for any map $f: G \rightarrow U(1)$.  Note that the phase~(\ref{vafa1loop})
is invariant under coboundaries, {\it i.e.}, it descends to a well-defined
map on group cohomology.)

To recap, discrete torsion is a degree of freedom in string orbifolds,
measured by the group cohomology group $H^2(G, U(1))$,
that corresponds to weighting twisted sector contributions to orbifold
partition functions by phases.  

Historically discrete torsion has been extremely mysterious.
It does not have an immediately obvious explanation, and so for a time
it was viewed as something intrinsic to string theory or conformal field 
theory, a smoking gun for string theory distinguishing it from other 
possible theories of quantum gravity.

Since the original paper \cite{vafa1}, there have been many papers
written on discrete torsion.  Rather than try to describe all of
them, we shall only describe two followups that the physics
community has deemed
particularly important:
\begin{enumerate}
\item In \cite{vafaed}, C.~Vafa and E.~Witten argued that turning on
discrete torsion could obstruct supersymmetric moduli
(moduli often naively identified with Calabi-Yau moduli).
\item In \cite{doug1,doug2}, M.~Douglas argued that turning on discrete
torsion had the effect of projectivizing
equivariant structures on
D-brane worldvolumes ({\it i.e.}, projectivized equivariant
K-theory), and these projectivized equivariant structures
were related to the usual closed-string description of discrete torsion.
\end{enumerate}
Also, there is a general belief that discrete torsion is intimately 
connected with the $B$ field, a (local) two-form tensor potential with
a gauge symmetry closely analogous to that of connections on 
principal $U(1)$ bundles, namely
$B \mapsto B + d \Lambda$ for any one-form $\Lambda$ is a symmetry
of the theory.  (More generally, a $B$ field is a higher-tensor analogue
of a connection on a principal $U(1)$ bundle.
Phrased yet more formally, the $B$ field is a connection on a gerbe,
and readers are referred to \cite{hitchin} for a more thorough description.)
The precise relationship between the discrete torsion and the $B$ field
has been somewhat elusive in the past; however, any serious attempt
to understand discrete torsion is certainly expected to explain
the precise nature of this relationship, and whether  
the $B$ field itself is sufficient, or whether some 
conformal-field-theory-specific effects also play a role.

In this section we shall outline
a purely mathematical description of discrete torsion.
Specifically, we shall argue that:
\begin{quotation}
{\it Discrete torsion is the choice of orbifold group action on the
$B$ field.}
\end{quotation}
In other words, the $B$ field itself is sufficient to explain
discrete torsion, one need not invoke any conformal-field-theory-specific
effects, and more generally, string theory need not enter
the discussion in any meaningful way.
Technically, by considering orbifold group actions on $B$ fields one
recovers not only discrete torsion, but also some other, more obscure
degrees of freedom also associated with the $B$ field, but we shall 
concentrate on explaining discrete torsion.

In general, whenever you have a field with a gauge symmetry, specifying
the group action on the underlying space does {\it not} suffice to specify the
group action on the theory.  After all, one can combine the group
action with a gauge transformation.
So, you {\it must} specify the orbifold group action on the fields,
not just the space.  More formally, a choice of orbifold group action
is known as an equivariant structure, so one must pick equivariant
structures on all fields with gauge symmetries.

In the context of heterotic toroidal orbifolds, the choice of
orbifold group action on the gauge fields is often called
``orbifold Wilson lines.''  Similarly, we shall
outline how the choice of orbifold group action
on the $B$ field is what is known as discrete torsion.

\subsection{Orbifold group actions on principal $U(1)$ bundles with connection}

The
choice of orbifold group action on $U(1)$ gauge fields 
({\it i.e.}, principal $U(1)$ bundles with connection)
forms a precise prototype of discrete torsion,
and is well-understood within the mathematics community
(see \cite{kostant} for an early reference).
For example, we shall review below how
such orbifold group actions are classified by
$H^1(G, U(1))$, whereas discrete torsion is classified by
$H^2(G, U(1))$.
As the technical details of equivariant structures 
(orbifold group actions) on $U(1)$ gauge fields are
both closely analogous to and much simpler than those for discrete
torsion, we shall review this analogue before proceeding
to discrete torsion.

How can one see the $H^1(G, U(1))$ advertised?
First, let us review equivariant structures on principal bundles
with connection.
The elegant way to proceed is as follows.  Let $L$ denote a 
principal $U(1)$ bundle over a space $X$, then a choice of
equivariant structure (orbifold group action) on $L$ is
a lift of the action of $G$ to $L$, {\it i.e.}, for each
$g \in G$, one defines a map $g': L \rightarrow L$
making the following diagram commute:
\begin{equation}
\xymatrix{
L \ar[r]^{ g' } \ar[d] & L \ar[d] \\
X \ar[r]^{ g } & X
}
\end{equation}
and obeying the group law, {\it i.e.}, $(g_1 g_2)' = 
g_1' g_2'$.
An equivariant structure on a principal $U(1)$ bundle with connection
is defined with the added constraint that each lift $g'$ must
preserve the connection.  In general, it is well-known
that such equivariant structures
need not exist, and even when they do exist, they are
not unique.  We shall only be concerned with non-uniqueness
here.

Now, it is straightforward to see that the set of
equivariant structures on a principal $U(1)$ bundle with
connection is a torsor under $H^1(G, U(1))$.
Given two lifts $g'$, $g''$ of a fixed $g \in G$, 
$\phi_g \equiv g' \circ (g'')^{-1}$ is a base-preserving bundle automorphism,
{\it i.e.}, a gauge transformation, and the constraint that
each lift preserve the connection becomes the constraint
that the gauge transformation $\phi_g$ 
preserve the connection -- so if $X$ is connected,
$\phi_g$ is a constant map into $U(1)$.
Finally, the constraint that the lifts respect the group law
becomes the statement that $\phi_g$ respects\footnote{
Note that we have used the fact that the structure group
of the principal bundles is abelian in this step -- the
analogous statement is not true in general of, say,
principal $SU(2)$ bundles.} the group law.
Hence, the $\phi_g$ define an element of
$\mbox{Hom}(G, U(1)) = H^1(G, U(1))$.

Thus, as advertised we see that the difference between
any two equivariant structures on a principal $U(1)$ bundle
with connection is given by an element of $H^1(G, U(1))$,
and in fact it is easy to check that the set of
equivariant structures is a torsor under $H^1(G, U(1))$.
An analogous analysis of $B$ fields will yield $H^2(G, U(1))$
(together with some other degrees of freedom related to the $B$ field).

For later use in describing $B$ fields, we can also repeat
this analysis more mechanically, in terms of data assigned
to an open cover of $X$.
Let $\{ U_{\alpha} \}$ denote a good open cover of our manifold
$X$, i.e., each $U_{\alpha}$ looks like a (contractible) open ball
inside $X$.  Then, to each $U_{\alpha}$ we associate a one-form
$A^{\alpha}$.  The one-forms on overlapping $U_{\alpha}$'s must be
related by a gauge transformation, i.e.,
\begin{displaymath}
A^{\alpha} \: - \: A^{\beta} \: = \: d \log g_{\alpha \beta}
\end{displaymath}
for some gauge transformations $g_{\alpha \beta}:
U_{\alpha} \cap U_{\beta} \rightarrow U(1)$, which close on
triple overlaps as
\begin{displaymath}
g_{\alpha \beta} \, g_{\beta \gamma} \, g_{\gamma \alpha}
\: = \: 1
\end{displaymath}

Next, we need to define the action of $G$ on such a $U(1)$ gauge
field.  More precisely, we need to relate\footnote{Technically,
for simplicity of presentation
we are assuming the elements of cover $\{ U_{\alpha} \}$ are invariant
under the action of $G$.  So, $\{ U_{\alpha} \}$ is not a good cover,
but we can assume each $U_{\alpha}$ is a disjoint union of
contractible open sets, which amounts to the next best thing.} $g^* A^{\alpha}$
to $A^{\alpha}$ and $g^* g_{\alpha \beta}$ to $g_{\alpha \beta}$.
The form of such orbifold group actions is a standard result:
\begin{eqnarray*}
g^* A^{\alpha} & = & A^{\alpha} \: + \: d \log h^g_{\alpha} \\
g^* g_{\alpha \beta} & = & \left(g_{\alpha \beta}\right) \, 
\left( h^g_{\alpha} \right) \, \left( h^g_{\beta} \right)^{-1} \\
h^{g_1 g_2}_{\alpha} & = & \left( g_2^* h^{g_1}_{\alpha} \right) \,
\left( h^{g_2}_{\alpha} \right)
\end{eqnarray*}
for each $g \in G$, for some collection of maps $h^g_{\alpha}: U_{\alpha}
\rightarrow U(1)$, which define the specific $G$-action.
Note, for example, that the first line is merely the statement
that $g^* A^{\alpha}$ and $A^{\alpha}$ need only agree up to a gauge
transformation; in order to specify the $G$-action, one must specify
those gauge transformations.

Finally, in order to see $H^1(G, U(1))$ explicitly, we need to study
the differences between distinct $G$-actions.
Let $\left( h^g_{\alpha} \right)$ define one $G$-action,
and $\left( \overline{h}^g_{\alpha} \right)$ another, on the same
$U(1)$ gauge field.
Define
\begin{displaymath}
\phi^g_{\alpha} \: = \: \frac{ h^g_{\alpha} }{
\overline{h}^g_{\alpha} }
\end{displaymath}
so that the $\phi^g_{\alpha}$ literally define the difference between
$G$-actions.  It is straightforward to check that the $\phi^g_{\alpha}$
determine a group homomorphism $G \rightarrow U(1)$, as we
shall outline below:
\begin{enumerate}
\item 
First use the fact that
\begin{eqnarray*}
g^* g_{\alpha \beta} & = &
 \left(g_{\alpha \beta}\right) \,
\left( h^g_{\alpha} \right) \, \left( h^g_{\beta} \right)^{-1} \\
\mbox{also } & = &  \left(g_{\alpha \beta}\right) \,
\left( \overline{h}^g_{\alpha} \right) \, \left( 
\overline{h}^g_{\beta} \right)^{-1} 
\end{eqnarray*}
By dividing these two lines, we see that $\phi^g_{\alpha} = 
\phi^g_{\beta}$ on $U_{\alpha} \cap U_{\beta}$, hence for any fixed
$g$ the $\phi^g_{\alpha}$ define a global function we shall call
$\phi^g$.
\item 
Next, use 
\begin{eqnarray*}
g^* A^{\alpha} & = & A^{\alpha} \: + \: d \log h^g_{\alpha} \\
\mbox{also } & = & A^{\alpha} \: + \: d \log
\overline{h}^g_{\alpha}
\end{eqnarray*}
to see that $\phi^g$ must be a constant function.
\item 
Finally, use
\begin{eqnarray*}
h^{g_1 g_2}_{\alpha} & = & \left( g_2^* h^{g_1}_{\alpha} \right) \,
\left( h^{g_2}_{\alpha} \right) \\
\overline{h}^{g_1 g_2}_{\alpha} & = & \left( g_2^* 
\overline{h}^{g_1}_{\alpha} \right) \,
\left( \overline{h}^{g_2}_{\alpha} \right)
\end{eqnarray*}
to see that $\phi^{g_1 g_2} = \phi^{g_1} \phi^{g_2}$.
\end{enumerate}
Thus, the $\phi^g$ define a group homomorphism $G \rightarrow U(1)$,
i.e., $\phi^g \in \mbox{Hom}(G, U(1))$.
However, it is a true fact that 
\begin{displaymath}
H^1(G, U(1)) \: = \: \mbox{Hom}(G, U(1))
\end{displaymath}
so again we see that the difference between any two orbifold group actions
on a $U(1)$ gauge field is defined by an element of $H^1(G, U(1))$.

In passing, we should mention that, if one is only interested
in finding group cohomology, there is a faster way to get it,
by using the Cartan-Leray spectral sequence description of
$G$-equivariant cohomology $H^2_G({\bf Z})$.  Unfortunately,
this cohomology class is not precisely directly physically relevant --
this gives information concerning all equivariantizable bundles
with all possible equivariant structures, yet we are concerned
with the set of equivariant structures on a {\it fixed} 
equivariantizable bundle {\it with connection}.  Also, using
the Cartan-Leray spectral sequence in this context 
obscures some important information; for example, the fact that
not all bundles are equivariantizable is hidden.
Partly the issue here is one of mere style;
however, considering that later we shall want to find more than
just group cohomology, we shall not quote the Cartan-Leray spectral
sequence when discussing discrete torsion.

\subsection{$B$ fields and $H^2(G, U(1))$}

When describing orbifold group actions on $B$ fields, one finds
$H^2(G, U(1))$, just as we found $H^1(G, U(1))$ for $U(1)$ gauge fields,
together with some related degrees of freedom.

As for bundles with connection, there is an elegant formal description of 
equivariant structures on gerbes with connection ({\it i.e.}, $B$ fields).
Let $P$ denote a gerbe over a space $X$.  A lift of $g \in G$ from $X$
to $P$ is given by a map $g': P \rightarrow P$ making the following
diagram commute:
\begin{equation}
\xymatrix{
P \ar[r]^{g'} \ar[d] & P \ar[d] \\
X \ar[r]^{ g } & X
}
\end{equation}
As before, the lifts $g'$ must obey the group law,
meaning that
\begin{equation}
\xymatrix{
P \ar[r]^{ g_1' } \ar@/_1pc/[rr]_{ (g_1 g_2)' } &
P \ar[r]^{ g_2' } & P
}
\end{equation}
Now, there is an additional layer of subtlety.  
If we realize the gerbe $P$ as a stack, then the diagram
above is a commutative diagram of (sheaves of) functors.
For a diagram of functors to commute does not mean that the
compositions must be equal, but merely isomorphic.
Hence, we must also specify isomorphisms
$\omega(g_1, g_2): g_2' \circ g_1' \Longrightarrow (g_1 g_2)'$,
and to be consistent on triples, we must demand
\begin{displaymath}
\omega(g_1 g_2, g_3) \circ \omega(g_1, g_2) \: = \:
\omega(g_1, g_2 g_3) \circ \omega(g_2, g_3)
\end{displaymath}
As before, an equivariant structure on a gerbe with connection
is defined with added constraints that the connection must
be preserved by these maps.
As for bundles with connection, $H^2(G, U(1))$ (and related
degrees of freedom) emerge when considering differences between
orbifold group actions,

As the details of this more elegant description are unfortunately
rather lengthy to work out, we shall instead proceed immediately
to a description of the $B$ field in terms of data assigned
to elements of an open cover.
Let $\{ U_{\alpha} \}$ be a cover of our manifold $X$ as before, 
then globally the $B$ field is described by \cite{hitchin}
a collection of two-forms $B^{\alpha}$ assigned to patches $U_{\alpha}$,
one-forms $A^{\alpha \beta}$ on double overlaps,
and $U(1)$-valued functions $h_{\alpha \beta \gamma}$ on
triple overlaps, satisfying
\begin{eqnarray*}
B^{\alpha} \: - \: B^{\beta} & = & d A^{\alpha \beta}  \\
A^{\alpha \beta} \: + \: A^{\beta \gamma} \: + \:
A^{\gamma \alpha} & = & d \log h_{\alpha \beta \gamma} \\
\left( h_{\alpha \beta \gamma} \right) \,
\left( h_{\alpha \beta \delta} \right)^{-1} \,
\left( h_{\alpha \gamma \delta} \right) \,
\left( h_{\beta \gamma \delta} \right)^{-1} & = & 1
\end{eqnarray*}
For physicists reading this discussion,
we should stress that the data $A^{\alpha \beta}$ and $h_{\alpha
\beta \gamma}$ are {\it not} some new fields in the theory,
just as bundle transition functions are not new fields in
a gauge theory;
rather, we are merely making gauge transformations on overlaps explicit.

Next, we need to define the $G$-action completely, which is to say,
describe how it acts not only on the $B^{\alpha}$ but also on
the overlap data $A^{\alpha \beta}$ and $h_{\alpha \beta \gamma}$.
The result can be derived from self-consistency, and also
exists in the literature (see for example \cite{brylinski}):
\begin{eqnarray*}
g^* B^{\alpha} & = & B^{\alpha} \: + \: d \chi(g)^{\alpha} \\
g^* A^{\alpha \beta} & = & A^{\alpha \beta} \: + \:
d \log \nu^g_{\alpha \beta} \: + \: \chi(g)^{\alpha} \: - \:
\chi(g)^{\beta} \\
g^* h_{\alpha \beta \gamma} & = & 
\left( h_{\alpha \beta \gamma} \right) \, 
\left( \nu^g_{\alpha \beta} \right) \,
\left( \nu^g_{\beta \gamma} \right) \,
\left( \nu^g_{\gamma \alpha} \right) \\
\chi(g_1 g_2)^{\alpha} & = &
\chi(g_2)^{\alpha} \: + \: g_2^* \chi(g_1)^{\alpha} \: - \:
d \log h^{g_1, g_2}_{\alpha} \\
\nu^{g_1 g_2}_{\alpha \beta} & = &
\left( \nu^{g_2}_{\alpha \beta} \right) \,
\left( g_2^* \nu^{g_1}_{\alpha \beta} \right) \,
\left( h^{g_1, g_2}_{\alpha} \right) \,
\left( g^{g_1, g_2}_{\beta} \right)^{-1} \\
\left( h^{g_1, g_2 g_3}_{\alpha} \right) \,
\left( h^{g_2, g_3}_{\alpha} \right) & = &
\left( g_3^* h^{g_1, g_2}_{\alpha} \right) \,
\left( h^{g_1 g_2, g_3}_{\alpha} \right)
\end{eqnarray*}
where $\chi(g)^{\alpha}$ (local one-forms),
$\nu^g_{\alpha \beta}$ (maps $U_{\alpha} \cap U_{\beta} \rightarrow U(1)$),
and $h^{g_1, g_2}_{\alpha}$ (maps $U_{\alpha} \cap U_{\beta} \cap
U_{\gamma} \rightarrow U(1)$) define the $G$ action on the $B$ field.

This looks more complicated than the description of $G$-actions on
$U(1)$ gauge fields, but the basic idea is identical.
For example, note from the first line that $g^* B^{\alpha}$
and $B^{\alpha}$ differ by a gauge transformation (defined by
$\chi(g)^{\alpha}$), for any $g \in G$.  Since the overlap
data is more complicated, one has to work harder to express
the complete $G$-action, but the basic principle is the same.

Now that we have defined $G$-actions on $B$ fields, we can
discuss the differences between $G$-actions on a fixed $B$ field.
When we did this for $U(1)$ gauge fields, we found that possible
$G$-actions are classified by $H^1(G, U(1))$.
Here, we shall find $H^2(G, U(1))$ (among other things).

We shall only outline how this section of the analysis proceeds.
Define, for example,
\begin{displaymath}
T^g_{\alpha \beta} \: \equiv \: \frac{ \nu^g_{\alpha \beta}}{
\overline{\nu}^g_{\alpha \beta} }
\end{displaymath}
Then using the relations
\begin{eqnarray*}
g^* h_{\alpha \beta \gamma} & = &
\left( h_{\alpha \beta \gamma} \right) \,
\left( \nu^g_{\alpha \beta} \right) \,
\left( \nu^g_{\beta \gamma} \right) \,
\left( \nu^g_{\gamma \alpha} \right) \\
\mbox{also } & = & \left( h_{\alpha \beta \gamma} \right) \,
\left( \overline{\nu}^g_{\alpha \beta} \right) \,
\left( \overline{\nu}^g_{\beta \gamma} \right) \,
\left( \overline{\nu}^g_{\gamma \alpha} \right) 
\end{eqnarray*}
we see that 
\begin{displaymath}
T^g_{\alpha \beta} \, T^g_{\beta \gamma} \, T^g_{\gamma \alpha} \: = \: 1
\end{displaymath}
so the $T^g_{\alpha \beta}$ are transition functions for a bundle,
call it $T^g$, for each $g \in G$.  Similarly,
$\chi(g)^{\alpha} - \overline{\chi}(g)^{\alpha}$ defines a flat
connection (a flat $U(1)$ gauge field) on $T^g$, and 
$h^{g_1, g_2}_{\alpha} / \overline{h}^{g_1, g_2}_{\alpha}$
defines an isomorphism $T^{g_2} \otimes g_2^* T^{g_1} \rightarrow T^{g_1 g_2}$,
satisfying a consistency condition we shall describe momentarily.

A moment's reflection reveals that this story is closely analogous
to the case of $G$-actions on $U(1)$ gauge fields.
There, recall the difference between any two $G$-actions was defined
by a set of gauge transformations $\phi^g$, one for each $g \in G$,
respecting the group law in $G$.  Here we see the same thing.
After all, globally a gauge transformation of a $B$ field is defined
by a $U(1)$ gauge field, so again we see the difference between two
$G$-actions is given by a set of gauge transformations, here
determined by the $T^g$.  The isomorphisms $T^{g_2} \otimes g_2^* T^{g_1}
\stackrel{\sim}{\longrightarrow} T^{g_1 g_2}$
simply enforce the group law on these gauge transformations.

To summarize our results so far, the difference between two $G$-actions
on a $B$ field is defined by
\begin{enumerate}
\item Bundles $T^g$ with flat connection (i.e., flat $U(1)$ gauge fields)
\item Maps $\omega^{g,h}: T^h \otimes h^* T^g \stackrel{\sim}{\longrightarrow}
T^{gh}$ such that the following diagram commutes:
\begin{equation}   \label{omegaconsis}
\xymatrix{
T^{g_3} \otimes g_3^* \left( \, T^{g_2} \otimes g_2^* T^{g_1} \, \right)
\ar[d]_{ \omega^{g_2, g_3} } 
\ar[r]^-{\omega^{g_1, g_2} } &
T^{g_3} \otimes g_3^* T^{g_1 g_2} 
\ar[d]^{\omega^{g_1 g_2, g_3} } \\
T^{g_2 g_3} \otimes (g_2 g_3)^* T^{g_1} 
\ar[r]^-{\omega^{g_1, g_2 g_3}} & T^{g_1 g_2 g_3}
}
\end{equation}
\end{enumerate}

Now we can finally see $H^2(G, U(1))$.
To make this explicit, take all of the bundles $T^g$ to be canonically
trivial with zero connection (i.e., set all of the $U(1)$ gauge fields
to zero), then the maps $\omega^{g,h}$ are forced to be constant gauge
transformations.  Commutivity of diagram~(\ref{omegaconsis})
becomes precisely the group 2-cocycle condition~(\ref{2cocycle}).
There are also residual gauge transformations, namely
constant gauge transformations of the bundles $T^g$, and these
act merely to change the $\omega^{g,h}$ by coboundaries.

Of course, not all flat bundles $T^g$ are trivializable,
and not all flat connections on trivializable bundles are gauge trivial.
So, clearly there are extra degrees of freedom present besides
merely $H^2(G, U(1))$, and a more careful analysis reveals
these are precisely momentum-dependent phase factors,
a more obscure aspect of string orbifolds also associated with
$B$ fields.

\subsection{Vafa's phases}

So far we have derived the classifying group $H^2(G, U(1))$ of
discrete torsion from orbifold group actions on $B$ fields.
In this section we will explain how one can see
the orbifold partition function phases that we originally
used to motivate discrete torsion.

\begin{figure}
\centerline{\psfig{file=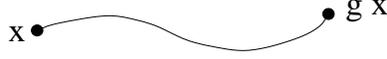,width=2in}}
\caption{ \label{figWline}
Lift of closed loop in $X/G$ to covering space $X$.
}
\end{figure}

These phases are closely analogous to orbifold Wilson lines, so again
let us briefly review relevant aspects of orbifold Wilson lines.
Consider a path on the covering space $X$, whose ends are related by
the action of $G$, as illustrated in figure~\ref{figWline}.
In other words, consider a point-particle one-loop
twisted sector, something which, on the quotient space $X/G$, 
becomes a closed loop.

If one computes the Wilson loop about that closed loop, but
upstairs on $X$, one finds it has the form
\begin{displaymath}
\varphi \, \exp \left( \, \int_x^{g \cdot x} A \, \right)
\end{displaymath}
Mostly this is the holonomy along the path in figure~\ref{figWline},
but there is a correcting factor $\varphi$ that implicitly describes
the $G$-action on $A$, i.e., $\varphi$ relates $A_x$ to $A_{g \cdot x}$.

Similar considerations for the $B$ field holonomy $\exp\left( B \right)$ on
polygons with sides identified by $G$ yield Vafa's twisted sector phases.

\begin{figure}   
\centerline{\psfig{file=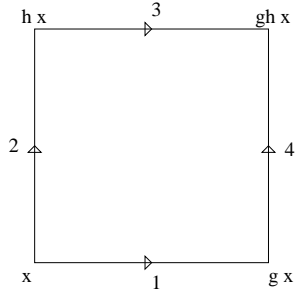,width=1.5in}}
\caption{  \label{fig1loop}
A twisted sector
contribution to the one-loop partition function.}
\end{figure}

Consider a one-loop twisted sector (i.e., from a string orbifold
on $T^2$) as shown in figure~\ref{fig1loop}.
Naively if we calculate the holonomy of the $B$ field about the
$T^2 \subset X/G$ one might expect
\begin{equation}
\exp \, \left( \, \int_{S} B \, \right)
\end{equation}
corresponding to the holonomy over the interior $S$ of the polygon.
However, this omits the contribution from gauge transformations at
the boundaries.

How can we take into account such gauge transformations?
Recall that a $G$ action on the $B$ field is specified by
\begin{itemize}
\item principal $U(1)$ bundles $T^g$ with flat connection $\Lambda(g)$, for all $g \in G$
\item connection-preserving bundle maps $\omega^{g,h}: T^h \otimes h^* T^g \stackrel{ \sim }{
\longrightarrow} T^{gh}$
enforcing the group law
\end{itemize}
Since we have gauge fields and the boundary of the square in 
figure~\ref{fig1loop} has one-dimensional components, the first thing
to try is to add a factor corresponding to the Wilson lines of the
$U(1)$ gauge fields along the boundaries, as
\begin{equation}    \label{try2}
\exp\left( \, \int_x^{h \cdot x} \Lambda(g) \: - \:
\int_x^{g \cdot x} \Lambda(h) \, \right) \,
\exp\left( \, \int_S B \, \right)
\end{equation}
(Signs are determined by relative orientations.)

However, expression~(\ref{try2}) is still not right, for technical reasons
(such as the fact that the result is not invariant under gauge
transformations of the $U(1)$ gauge fields $\Lambda$).
As described in much more detail in \cite{dt3}, to fix this expression
we must also take into account the corners, which yields the
correct general expression
\begin{equation}   \label{Bholon}
\left( \, \omega_x^{g,h} \, \right) \,
\left( \, \omega_x^{h,g} \, \right)^{-1} \,
\exp\left( \, \int_x^{h \cdot x} \Lambda(g) \: - \:
\int_x^{g \cdot x} \Lambda(h) \, \right) \,
\exp\left( \, \int_S B \, \right)
\end{equation}
(Overall normalizations are fixed by comparison to $B$ field holonomies
around nontrivial cycles on the quotient space.)

Now, for those degrees of freedom measured by $H^2(G, U(1))$,
recall we can assume $\Lambda(g) \equiv 0$ for all $g \in G$
and that $\omega^{g,h}$ is constant, hence for these degrees of freedom
the expression~(\ref{Bholon}) reduces to
\begin{equation}
\left( \, \omega^{g,h} \, \right) \,
\left( \, \omega^{h,g} \, \right)^{-1} \,
\exp\left( \, \int_S B \, \right)
\end{equation}
Note that the resulting phase factor completely agrees with Vafa's 
one-loop phase factor as stated in equation~(\ref{vafa1loop}).

Since all contributions to the path integral (in this twisted sector)
are weighted with this same phase, the effect is to multiply
the partition function for this twisted sector by Vafa's phase.
Hence, we find complete agreement with Vafa's original description
\cite{vafa1}.

Other checks of this description, such as multiloop
factorization, also work out nicely and are described in
\cite{dt3}.

\subsection{Douglas's projectivization for D-branes}

An action of $G$ on a $U(N)$ gauge field is described by
\begin{eqnarray}
g^* A^{\alpha} & = &
\left( \gamma^g_{\alpha} \right) \, A^{\alpha} \,
\left( \gamma^g_{\alpha} \right)^{-1} \: + \:
\left( \gamma^g_{\alpha} \right) \, 
d \left( \gamma^g_{\alpha} \right)^{-1} \\
g^* g_{\alpha \beta} & = &
\left( \gamma^g_{\alpha} \right) \,
\left( g_{\alpha \beta} \right) \,
\left( \gamma^g_{\beta} \right)^{-1} \\
\gamma^{gh}_{\alpha} & = & 
\left( h^* \gamma^g_{\alpha} \right) \,
\left( \gamma^h_{\alpha} \right)    \label{untw}
\end{eqnarray}
for some $\gamma^g_{\alpha}: U_{\alpha} \rightarrow U(N)$
defining the $G$-action.

M. Douglas conjectured \cite{doug1,doug2} that when discrete torsion is turned
on, the $G$-action on a D-brane gauge field is twisted,
meaning that equation~(\ref{untw}) is replaced with
\begin{equation}  \label{tw1}
\left( \omega^{g,h} \right) \,
\left( \gamma^{gh}_{\alpha} \right) \: = \:
\left( h^* \gamma^g_{\alpha} \right) \,
\left( \gamma^h_{\alpha} \right)  
\end{equation}
so that instead of an honest representation of the orbifold
group $G$, one only has a projective representation.
Phrased another fashion, instead of working with equivariant K-theory,
one works with projectivized equivariant K-theory.

This projectivized representation comes from the fact that in the
presence of a nontrivial $B$ field, the ``bundle'' on the worldvolume
of a D-brane is twisted\footnote{Technically, we should mention that
the equivariant structure described below is almost but
not quite uniquely fixed
by self-consistency; rather, this reflects some minor minimal choices,
and is the analogue of a `true' equivariant structure as opposed,
for example, to a projective equivariant structure.}:
\begin{eqnarray*}
A^{\alpha} \: - \:
g_{\alpha \beta} A^{\beta} g_{\alpha \beta}^{-1} 
\: - \: d \log g_{\alpha \beta}^{-1} & = &
A^{\alpha \beta} I \\
g_{\alpha \beta} \, g_{\beta \gamma} \, g_{\gamma \alpha}
& = &  h_{\alpha \beta \gamma} I
\end{eqnarray*}
where $A^{\alpha}$ is a local $U(N)$ gauge field and
$\left( B^{\alpha}, A^{\alpha \beta}, h_{\alpha \beta \gamma}
\right)$ define the $B$ field (as discussed earlier).
This twisting is a consequence of the fact that under
$B \mapsto B + d \Lambda$, the Chan-Paton gauge field
$A \mapsto A - \Lambda I$, as described in for example
\cite{freeded}.

As a direct result of this intermingling between the gauge field
and the $B$ field, any $G$-action on the ``bundle'' must be
intertwined with the $G$-action on the $B$ field.

Using self-consistency arguments, we find that in general,
the $G$-action on the D-brane ``bundle'' is actually of the
form
\begin{eqnarray}
g^* A^{\alpha} & = &
\left( \gamma^g_{\alpha} \right) \, A^{\alpha} \,
\left( \gamma^g_{\alpha} \right)^{-1} \: + \:
\left( \gamma^g_{\alpha} \right) \,
d \left( \gamma^g_{\alpha} \right)^{-1} 
\: + \: \chi(g)^{\alpha} I \\
g^* g_{\alpha \beta} & = &
\left( \nu^g_{\alpha \beta} \right) \,
\left[ \,
\left( \gamma^g_{\alpha} \right) \,
\left( g_{\alpha \beta} \right) \,
\left( \gamma^g_{\beta} \right)^{-1} \, \right] \\
\left( \omega^{g,h}_{\alpha} \right) \, 
\left( \gamma^{gh}_{\alpha} \right) & = &
\left( h^* \gamma^g_{\alpha} \right) \,
\left( \gamma^h_{\alpha} \right)  \label{tw2}
\end{eqnarray}
where $\chi(g)^{\alpha}$, $\nu^g_{\alpha \beta}$, and $\omega^{g,h}_{\alpha}$
are data defining the $G$-action on the $B$ field.

To be very brief, from equation~(\ref{tw2}) we see Douglas's 
projectivization (suitably generalized).

\subsection{Discrete torsion for 3-form potentials}

So far we have only discussed $B$ fields.
However, string theory and $M$ theory have other tensor field potentials.
These other tensor field potentials have precise analogues of 
discrete torsion \cite{cdt}.

For example, consider a three-form potential $C_{\mu \nu \rho}$. 
It can be shown \cite{cdt} that just as $G$-actions on $U(1)$
gauge fields are classified by $H^1(G, U(1))$, and $G$-actions
on $B$ fields are (partially) classified by $H^2(G, U(1))$,
possible $G$-actions on $C$ fields are (partially)
classified by $H^3(G, U(1))$.

\begin{figure}
\centerline{\psfig{file=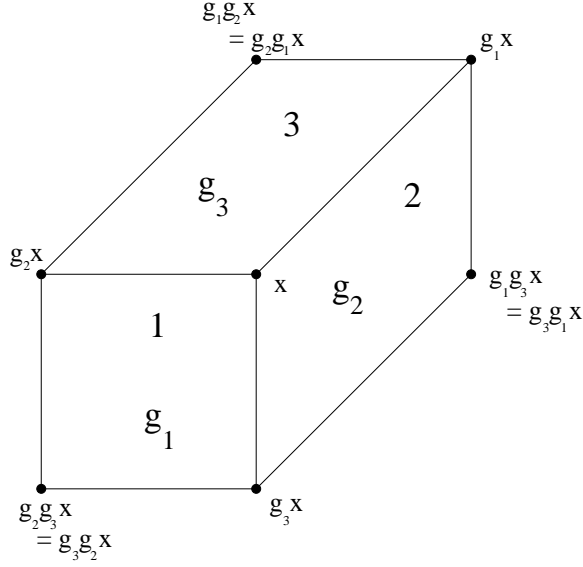,width=3in}}
\caption{ \label{figcbox}
Three-torus seen as open box on covering space.}
\end{figure}

One can also get phase factors for membranes, closely analogous
to Vafa's phase factors.
For example, recall that for a $T^2$ twisted sector
(as illustrated in figure~\ref{fig1loop}), from $\exp \left(
\int B \right)$ one has a phase factor
\begin{equation}
\omega(g,h) \: - \: \omega(h,g)    \label{vafaorig}
\end{equation}
(where we have chosen here to write the abelian product additively instead
of multiplicatively).
In the present case, for membranes on $T^3$ twisted sectors as illustrated
in figure~\ref{figcbox}, from $\exp \left( \int C \right)$ one
has a phase factor
\begin{equation}
\omega(g_1, g_2, g_3) \: - \: \omega(g_2, g_1, g_3) \: - \:
\omega(g_3, g_2, g_1) \: + \: \omega(g_3, g_1, g_2) \: + \:
\omega(g_2, g_3, g_1) \: - \: \omega(g_1, g_3, g_2)
\label{memphase}
\end{equation}

Just as Vafa's original phase factor~(\ref{vafaorig})
for $T^2$ was $SL(2, {\bf Z})$ invariant, the phase 
factor~(\ref{memphase}) for $T^3$ is $SL(3, {\bf Z})$ invariant
(reflecting the fact that it is well-defined on $T^3$).

\subsection{Vafa-Witten and supersymmetric moduli}

In the paper \cite{vafaed} it was pointed out that
in string orbifolds with discrete torsion `turned on,'
supersymmetric moduli are often obstructed.
It should be mentioned that because in many simple cases,
supersymmetric moduli can be identified with Calabi-Yau moduli,
the paper \cite{vafaed} has sometimes been misinterpreted to mean that
Calabi-Yau are somehow obstructed, that the geometry of the
Calabi-Yau somehow changes, but in fact the authors of \cite{vafaed} 
argued only 
a much weaker statement.

One of their proposed explanations for this behavior now
seems extremely reasonable (for reasons explained in much
greater detail in \cite{dt3}).  Namely, the authors of
\cite{vafaed} speculated that,
after `turning on' discrete torsion, attempts to resolve or deform
singularities will often result in nonzero $B$ field curvature
(known to physicists as `torsion,' and denoted by $H$).
(Such behavior is closely analogous to one description of
the McKay correspondence \cite{gsv}, for example.)
Nonzero $B$ field curvature breaks supersymmetry, and so such
deformations would be obstructed.

Although this sounds extremely plausible in the present context,
it has not yet been verified to date, and is essentially the only
remaining aspect of classical discrete torsion that is not yet
completely understood.
(See however \cite{berenleigh} for a recent attempt to 
understand this effect on D-branes using noncommutative geometry.)

\section{Quotient stacks and string orbifolds}

In this first half of this paper we gave a largely complete
understanding of discrete torsion, as the choice of
equivariant structure on the $B$ field (a gerbe with connection).
In the second half, we shall shift gears and describe the
relation between string orbifolds and quotient stacks.
Although in the first half of this paper we were able to
give a fairly complete understanding of discrete torsion, in the second half
we shall only set up the basics required to relate string orbifolds
to quotient stacks, and emphasize that much work remains to be done.

As mentioned in the introduction, historically physicists
have assumed that string orbifolds describe strings on
quotient spaces with some sort of `stringy'
behavior at singularities.  Of course, string orbifolds are not sigma models
on quotient spaces; rather, they describe group actions on covers.
However, physicists have assumed that this was merely 
scaffolding.
In fact, in practice, although physicists often speak of `stringy'
behavior at singularities, when relating string orbifolds
to the rest of string theory, we often implicitly identify
string orbifolds with quotient spaces.
For example, descriptions of string moduli spaces, fundamental to
topics ranging from mirror symmetry to string/string duality,
all implicitly assume that the deformation theory of a string
orbifold is the same as that of a quotient space.

Unfortunately, thinking about string orbifolds in terms
of quotient spaces with `stringy' behavior at singularities
is not wholly satisfactory -- various physical features
of string orbifold CFT's, such as well-behavedness of the CFT,
have only clumsy explanations within the present physics lore. 
For such reasons, a more subtle alternative was proposed in \cite{qstx},
namely that string orbifolds do not describe strings on quotient spaces,
with or without `stringy' behavior,
but rather describe strings compactified on quotient {\it stacks}. 

To most physicists, the idea that string orbifold CFT's coincide with
CFT's for 
strings compactified on quotient stacks, a formal geometric structure
assigned to the group-actions-on-covers scaffolding,
is somewhat radical.  Certainly nothing of the sort has previously been
believed within the physics community.
Reference \cite{qstx} concentrated on trying to make the notion
palatable to physicists.  In a nutshell, such a description has nontrivial
physical implications.  For example, by thinking about
string orbifolds as sigma models on stacks, one can immediately
resolve a number of puzzling issues about the physics of string
orbifolds, which were not satisfactorily resolved within 
the physics lore.  Features often ascribed to `stringy' effects
can be seen to be easy consequences of geometrical features
of stacks.  This also has nontrivial consequences for the understanding
of string moduli spaces.

To mathematicians familiar with stacks, on the other hand,
the idea that string orbifold CFT's coincide with CFT's
for strings compactified on quotient stacks
is much more natural.  After all, one way of thinking about quotient
stacks is as an overcomplicated way to describe group actions on covers,
which is certainly the language that string orbifolds are phrased in.
Certainly one can efficiently manipulate group actions on covers
by working with quotient stacks.
Moreover, quotient stacks are more than just an overcomplicated
way of describing group actions on covers -- they also can
be interpreted as `generalized spaces.'  For example, one sometimes
hears\footnote{Though good references are unfortunately very difficult
to find.  For this reason, reference \cite{qstx} includes a lengthy
discussion of differential geometry on stacks in general,
and quotient stacks in particular.}
that it is possible to do differential geometry on quotient stacks.
So, someone who was not acquainted with the physics literature
might be led to (very naively) believe that the notion of string
compactification on a stack is sensible, and furthermore that a 
string orbifold coincides with a string on a quotient stack.

Unfortunately, the notion that that extra structure
of a `generalized space' has physical content, in the fashion above,
has not been justified.
Any competent physicist would point out that, not only
does this appear to contradict the standard physics lore, but also a tremendous
amount of work must be done to even check whether the idea of
string compactification on a stack is sensible,
much less reconcile it with the role string orbifolds play in the
rest of string theory, or resolve the basic physical contradictions that
appear.

For example, before one can make sense out of the statement that
string orbifolds are the same thing as
strings compactified on {\it quotient} stacks,
one must first describe string compactification on stacks in general,
and check whether this is even a sensible notion.
The fact that one can do differential geometry on stacks is a necessary
condition, but by no means is it sufficient.
One way to do this (which we shall outline the basics of)
is to write down the classical action for a sigma model on a stack,
and then try to check whether it can be quantized (which involves
studying global issues).  As the second half of such requirements is
rather technical, one uses consistency checks to gain insight into
whether or not this is reasonable.
At the end of the day, in order to claim that one completely understands
these matters, one must be able to answer questions ranging from basics such as
``what is a string on another stack,
{\it e.g.} a gerbe''
to more difficult ones such as ``how does one make sense of sigma model
anomalies in this context'' ({\it i.e.}, ``what about global issues'')
and ``sigma models with stack targets are one thing,
but can do one quantum field theory directly on a stack, now viewed as
spacetime itself?''  ``What is the propagator for a scalar field
on $[R^4/Z_2]$?'' and so forth.

Even after one attempts to make sense of the notion of a string on a stack,
and checks whether a string orbifold really is a sigma model on 
a quotient stack, one still must reconcile quotient stacks with
the role string orbifolds play in the rest of string theory.
The statement that string orbifold CFT's coincide with CFT's for
strings compactified on quotient stacks has nontrivial physical
implications, which cannot be ignored.
For example, as mentioned above, in constructing moduli spaces (important
to justify everything from mirror symmetry to string/string duality),
physicists have assumed that the deformation theory of
string orbifolds was that of a quotient space,
as indeed we assumed string orbifolds described strings on quotient spaces.
If string orbifolds are indeed sigma models on quotient stacks,
then one must explain how the deformation theory agrees,
despite the fact that deformation theory of a quotient stack
is not the same as that of a quotient space \cite{pantevpriv}.

Perhaps the best question one can ask is simply:  why bother?
If, at the end of the day, describing string orbifolds in terms of
quotient stacks accomplishes nothing more than providing an
overcomplicated description of group actions on covers, then there
is hardly a point.  However, we shall argue that such a statement
has nontrivial physical implications.  First, by thinking about
string orbifolds as sigma models, one gains a much clearer understanding
of certain physical features of string orbifolds.  Second,
as mentioned above, such a statement has nontrivial implications for our
understanding of string moduli spaces -- if a string orbifold really is
a sigma model on a quotient stack (assuming that is a sensible notion),
then any deformations of the CFT must be understood in terms of
deformations of the quotient stack, not the quotient space.
 
In the next few subsections we shall begin by describing classical sigma
models on stacks, verify that string orbifolds are indeed
sigma models on quotient stacks, and describe strings on some
other stacks.  We shall discuss the massless spectrum of sigma models
on stacks, and in particular resolve the apparent contradiction that
the massless spectrum is not given by the cohomology of the target.
We shall also outline how one can understand the well-behavedness
of string orbifold CFT's in terms of the geometry of stacks.
Finally we shall conclude with some comments on and questions
about deformation theory.

More generally there seems to be a considerable gap between
the mathematics lore and the physics lore on string orbifolds.
We shall address these distinctions as they arise,
in an attempt to help bridge certain gaps.


Physicists who are not comfortable with quotient stacks
are encouraged to read \cite{qstx}, where we have spoken to
their concerns in print.  As a description for physicists
already exists, the rest of this
lecture shall instead be oriented towards mathematicians,
and in particular, that aspect of the mathematics community
which feels that they already know that string orbifolds
are the same as strings compactified on quotient stacks.

We should emphasize that we do not wish to claim that this
matter is completely resolved -- there is still a tremendous
amount of work that must be done to verify that string orbifolds
really can be consistently interpreted as sigma models on
quotient stacks, to verify that this is not only itself reasonable
but consistent with the role string orbifolds play in
string theory as a whole.  Although we shall outline some of the
basic work required to make sense of such notions, and resolve some of the
basic paradoxes that crop up, there are still some
strong arguments that ultimately
this program must fail, which we have not yet been able to resolve.
In other words, at present, despite what some might like to believe, 
there is still no guarantee that string orbifolds really do describe
strings on quotient stacks.

\subsection{Sigma models on stacks}

As described above, before one can say that string orbifolds
are the same as strings compactified on quotient stacks, one must first 
describe string compactification on general stacks.
A necessary condition for this is the ability to do differential
geometry on stacks -- something one sometimes hears mentioned.
Reference \cite{qstx} contains a lengthy discussion of this matter,
something that seems to be largely missing from the literature.

Now, being able to do differential geometry on stacks is not nearly 
sufficient to
enable one to speak of describing strings on stacks.
One must be able to make sense out of the notion of a sigma model
on a stack, and answer a long list of related questions.
We shall outline classical actions for sigma models on stacks (see
\cite{qstx} for more details), answer some of the basic
questions one can ask \cite{qstx}, and list more questions that must be answered
before compactification on stacks can be described as well understood.

Let ${\cal F}$ be a stack, with atlas $X$.
(We shall only attempt to describe sigma models on stacks with
atlases.)  For readers not well-acquainted with stacks,
for $X$ to be an atlas for ${\cal F}$ implies that
\begin{itemize}
\item implicitly there is also a fixed map $X \rightarrow {\cal F}$
(which is required to be a surjective local homeomorphism)
\item for any space $Y$ and map $Y \rightarrow {\cal F}$,
the fibered product $Y \times_{ {\cal F} } X$ is an honest space, not a stack.
\end{itemize}

For example, if ${\cal F}$ is a space, not just a stack
(spaces are special cases of stacks), then ${\cal F}$ is its own atlas,
and $Y \times_{ {\cal F} } X = Y$ for any $Y$.
For another example, suppose ${\cal F}$ is a quotient stack $[X/G]$,
with $G$ discrete and acting by diffeomorphisms on a smooth space $X$.
In such a case, $X$ is an atlas for $[X/G]$.  In this case,
$Y \times_{ [X/G] } X$ is a principal $G$-bundle over $Y$,
partially specifying the map $Y \rightarrow [X/G]$,
and the projection map $Y \times_{ [X/G] } X \longrightarrow X$
is the $G$-equivariant map from the total space of the bundle to $X$,
specifying the rest of the map $Y \rightarrow [X/G]$.

Now, the natural description of
a sigma model with target ${\cal F}$, formulated on
(base) space $Y$, is a sum over equivalence classes\footnote{A sigma model
path integral is a sum over maps, after all, hence one must take
equivalence classes in order to make sense out of such a sum.} of maps
$Y \rightarrow {\cal F}$, weighted by
$\exp(iS)$, where the classical action $S$ is
formulated as follows.  Fix a map $\phi: Y \rightarrow
{\cal F}$.  If we let\footnote{
Note that since both $Y \times_{ {\cal F} } X$ and $X$ are ordinary
spaces, $\Phi$ is a map in the ordinary sense of the term.}
$\Phi: Y \times_{ {\cal F} } X \rightarrow X$ denote the
second projection map (implicitly encoding part of the map
$\phi: Y \rightarrow {\cal F}$), then the natural proposal for
the bosonic part of the classical
action for a sigma model on ${\cal F}$ is given by \cite{qstx}
\begin{equation}
\int d^2 \sigma \, \left( \pi_1^* \phi^* G_{\mu \nu} \right)
\pi_1^* h^{\alpha \beta} 
\left( \frac{ \partial \Phi^{\mu} }{ \partial \sigma^{\alpha} } \right)
\left( \frac{ \partial \Phi^{\nu} }{ \partial \sigma^{\beta} } \right)
\end{equation}
where $h^{\alpha \beta}$ is the worldsheet metric, $\phi^* G$ denotes
the pullback of the metric on ${\cal F}$ to $Y$ (metrics on ${\cal F}$
are described
in terms of their pullbacks), 
$\pi_1: Y \times_{ {\cal F} } X \rightarrow Y$ is the projection map,
and this action is integrated over
a lift\footnote{Sensible essentially because
the (projection) map $\pi_1: Y \times_{ {\cal F} } X \rightarrow Y$
is a surjective local homeomorphism.}
of $Y$ to $Y \times_{ {\cal F} } X$.

A few examples should help clarify this description:
\begin{enumerate}
\item Suppose ${\cal F}$ is an ordinary space.
Then the path integral is a sum over maps into that space,
and as $Y \times_{ {\cal F} } X = Y$ (taking the atlas $X$ to be
${\cal F}$ itself), we see that $\Phi = \phi$,
and so 
in this case the classical action proposed
above duplicates the usual classical action, as described
in equation~(\ref{stdclassact}), as well as the path integral sum.
Thus, the description above duplicates sigma models
on ordinary spaces.
\item Suppose ${\cal F} = [X/G]$, where $X$ is smooth and
$G$ is a nontrivial action of a discrete group by diffeomorphisms.
Then the path integral is a sum over equivalence classes of maps
$Y \rightarrow [X/G]$, which is to say, equivalence classes of
principal $G$-bundles on $Y$ together with $G$-equivariant maps
from the total space of the bundle into $X$.  It is easy to check
that the proposed classical action above duplicates the usual
classical action for a string orbifold.  Also, by summing over
(equivalence classes of) maps $Y \rightarrow {\cal F}$,
note we are summing over both twisted sectors as well
as maps within any given twisted sector.

Now, for each such map $\phi: Y \rightarrow {\cal F}$, 
there are $| G |$ lifts of $Y$ to the 
total space of the bundle (which is $Y \times_{ [X/G] } X$),
{\it i.e.}, $| G |$ twisted sector maps, as they usually appear in
physics.

Note we are only summing over equivalence classes of bundles,
not all possible twisted sector maps.  However, we can trivially
sum over all possible twisted sector maps, at the cost of
overcounting by $| G |$.  Hence, we can equivalently describe
this in terms of a sum over twisted sector maps, 
but weighted by $| G |^{-1}$.  Hence we recover both the path
integral sum and the overall multiplicative factor of $|G|^{-1}$
appearing in string orbifold partition functions ({\it e.g.}
equation~(\ref{zt2})).
\end{enumerate}

Thus, we see that the natural definition of a sigma model on a stack
duplicates not only sigma models on ordinary spaces,
but also string orbifolds when the target is a quotient stack, 
even down to the $| G |^{-1}$ factor
appearing in partition functions.
Also, note that the fact that the path integral sum duplicates both the
twisted sector sum and the functional integral within each
twisted sector is independent of our proposal for a classical action -- 
even if our proposed classical action is wrong, agreement between
path integral sums still holds true, and is a `smoking gun' for
some sort of interpretation as a sigma model, as emphasized in
\cite{qstx}.

We should take this opportunity to also note that this description
of sigma models does not make any assumptions concerning the dimension
of the base space $Y$ -- classically there are analogues of `string' orbifolds
in every dimension, all obtained precisely by gauging the action of
a discrete group on the target space of a sigma model.

So far we have only recovered known results; let us now try 
something new.  Suppose the target ${\cal F}$ is a gerbe.
For simplicity, we shall assume that ${\cal F}$ is the
canonical trivial $G$-gerbe on a space $X$.
Such a gerbe is described by the quotient stack $[X/G]$,
where the action of $G$ on $X$ is trivial.
Using the notion of sigma model on a stack as above,
one quickly finds that the path integral for this target space
is the same as the path integral for a sigma model on $X$,
up to an overall multiplicative factor (equal to the number of
equivalence classes of principal $G$-bundles on $Y$).
As overall factors are irrelevant in path integrals,
the result appears to be that a string on the canonical trivial
gerbe is the same as a string on the underlying space.
More generally, it is natural to conjecture that
strings on flat gerbes are equivalent to
strings on underlying spaces, but with flat $B$ fields.
In particular, such a result would nicely dovetail with the well-known
fact that
a coherent sheaf on a flat gerbe is equivalent to a `twisted' sheaf
on the underlying space, the same twisting that occurs in the presence
of a $B$ field.  (For physicists, this is an alternative to the
description in terms of modules over Azumaya algebras that has
recently been popularized \cite{kapustin}.)

So far we have only discussed classical actions for sigma models
on stacks, but there is much more that must be done before 
one can verify that the notion of a sigma model on a stack
is necessarily sensible.  In effect, we have only considered
local behavior, but in order to be sure this notion is sensible
after quantization, one also needs to consider global phenomena.
Such considerations were the source of much hand-wringing when
nonlinear sigma models on ordinary spaces were first introduced
(see for example \cite{mmn}), and must be repeated for stacks.
We have not performed such global analyses,
but instead shall perform several consistency checks.
For example, in the next section we will perform such a check
by examining the massless spectrum of the purported sigma model,
which ordinarily must coincide with some cohomology of the target space.
Interestingly enough, we will find that stacks seem to fail this
consistency check!  Although we will overcome this particular difficulty,
we will not have solutions to the puzzles posed by other failed
consistency checks we will describe later.

\subsection{Massless spectrum of a sigma model on a stack}


In \cite{qstx}, we argued that thinking about 
(supersymmetric) sigma models
on quotient stacks $[X/G]$ (with $G$ acting effectively)
led one to conclude that the massless spectrum
should be given by cohomology of the associated inertia group stack
$I_{[X/G]}$:
\begin{equation}   \label{oec}
I_{[X/\Gamma]} \: \cong \: \coprod_{[g]} \, \left[ \, X^g / C(g) \, \right]
\end{equation}
a result that, although obscure to most physicists,
is known to some mathematicians.
(Note the obvious relation to the Hirzebruch-H\"ofer \cite{hirzhofer}
description of
orbifold Euler characteristics, as in expression~(\ref{hhoec}).)

Note that we seem to immediately find a contradiction.
The massless spectrum of a string sigma model is given by some
cohomology of the target; yet, in a string orbifold,
the twisted sectors of the massless spectrum are not described
by cohomology of the (quotient stack) target.  On the face of it,
this completely contradicts the idea that string orbifolds can 
consistently be considered to be sigma models on quotient stacks,
and indeed, many physicists would take this as strong evidence
that string orbifolds {\it cannot} possibly be the same as strings
compactified on stacks.

In fact, a physics subtlety saves the day.
As observed in \cite{qstx}, the usual statement that the massless
spectrum of a string sigma model is the cohomology of the target
is a bit imprecise -- it would be better to say, the massless
spectrum of a string sigma model is the cohomology of the
zero-momentum part of the loop space of the target.
When the target is an ordinary space, the zero-momentum part of the
loop space is the original space itself.
However, when the target is a stack, there is a distinction.
Although we were able to resolve this apparent paradox,
we shall see more apparent contradictions later,
and unfortunately we will not be able to describe how to solve the
puzzles they pose.

Another important point to note is that we have been naturally
led to a description of twist fields that is very different
from the description most physicists assumed would hold true.
Indeed, the form of this
description of twist fields we have just given
is largely unknown
within the physics community, although after unraveling definitions,
it boils down to a description in terms of group actions on covers
that appears more familiar.
As described in section~\ref{lightning}, in the physics community
many assumed that twist fields could be understood in terms of
some cohomology of the quotient space, called `orbifold cohomology,'
that would implicitly encode information about the hypothesized
`stringy' behavior
at singularities.  Put another way, the physics of string orbifolds
seemed to predict the existence of such a description of twist fields,
not just in terms of group actions on covers, but in terms of 
some cohomology of a quotient space.
However, not only is the description above not a cohomology
of quotient spaces, it is not even a cohomology of a quotient stack!
By thinking about the physics of string orbifolds in terms of
stacks, a more subtle approach than most physicists have
previously considered, we have been led to a very different description of
twist fields than most physicists have assumed would hold true.
Put another way, by thinking in terms of stacks, one is led to
conclude that a description of twist fields in terms of group actions
on covers is the best one can hope to manage, that previous
expectations of some description directly in terms of the
underlying quotient space were naive.

Although expression~(\ref{oec}) is largely unknown within the
physics community (and indeed, substantially deviates from
what physicists have traditionally expected),
it is not unknown within the mathematics community.
For example, it is formally equivalent\footnote{
Technically the paper \cite{cr} used a description of
twist fields (in terms of group actions on covers)
equivalent to the one above as a starting point,
and concentrated on developing an ansatz for a cup product which,
by construction, duplicates the twist field correlation functions
studied in the older physics literature.  The distinction we are
trying to draw here appears to be somewhat more basic -- many physicists would
like some cohomology theory of a quotient space, from which twist fields
emerge naturally, whereas a description of twist fields in terms of
group actions on covers that begins with expression~(\ref{oec}) does  
not sound like a final answer to some physicists.} to the description of
twist fields given in \cite{cr}.  We should emphasize again,
however, that the definition of `orbifold cohomology'
presently used in the mathematics community appears to be somewhat
different
from the definition used by many members of the physics community.
A description of twist fields in terms of group actions on covers,
although interesting, appears neither to be the `orbifold cohomology'
that many physicists have traditionally desired, nor does it seem to shed
insight into the physics questions that physicists have
hoped orbifold cohomology would answer.
We have argued, from an understanding
of string orbifolds as sigma models on quotient stacks,
that \cite{cr} is nevertheless the correct and `final' understanding
of twist fields.  However, if future analysis reveals that we
are wrong, then although \cite{cr} is very interesting, 
some physicists would argue that it is 
not necessarily the final word on twist fields.

\subsection{Well-behavedness of string orbifold CFT's}

So far we have only set up {\it how} one could describe string orbifolds
in terms of stacks.  We have yet to explain {\it why} one would
wish to do so.  After all, if we are only using quotient stacks
as a radically overcomplicated way to describe group actions on covers,
then there is hardly a point.  However, this description has
nontrivial physical consequences, which is the real
reasons why physicists should be interested in such a description.
We shall describe in this section 
how this description greatly clarifies the physics of
string orbifolds.  Viewing string orbifolds
as sigma models on stacks sheds new light on
the physics of string orbifolds, and many properties
of string orbifolds that were previously attributed to
some sort of `stringy' effect can be seen
to be easy consequences of the geometry of stacks.  In later sections we
shall describe other nontrivial consequences of working with stacks.

For example,
historically one of the more interesting and useful
features of string orbifold
CFT's was their well-behavedness:  string
orbifolds were constructed in an attempt to describe
strings on quotient spaces, yet even when
the quotient space is singular, the CFT is well-behaved.
Historically physicists have often quoted two general mechanisms
in connection with this
behavior:
\begin{itemize}
\item ``String orbifold CFT's are well-behaved because a
CFT operator typically associated with holonomy of the $B$ field
on exceptional divisors (the `theta angle') is expected to often have
a nonzero vacuum expectation value''
\cite{edstrings95,paulz2,kw}.

Certainly the string orbifold
point in a moduli space of string vacua is distinct from the
point corresponding to a sigma model on the quotient space,
and one natural mechanism for this to occur is if the string
orbifold point corresponds to a nonzero theta angle.
Moreover, from a linear sigma model perspective this mechanism
is extremely natural \cite{edstrings95}.

Now, one might ask if this phenomenon has a more intrinsic
explanation.  This nonzero theta angle is described either
by studying sigma models on resolutions of the underlying quotient space,
and considering rational-curve-counting in limits when the exceptional
divisors shrink to zero size, or in terms of linear sigma models,
in which the physical theory is described very indirectly in terms
of a distinct theory which is believed to
flow (in the sense of the renormalization group) to the relevant
physical theory.  Neither of these descriptions involves the
string orbifold CFT directly, but rather only
deformations of that physical theory. 

If one wishes to propose a description of string orbifolds as 
something other than strings on quotient spaces with
stringy effects at singularities, then one natural question that
will be asked is, can one understand this nonzero theta angle
business geometrically?  We have described nonzero theta angles
in terms of CFT operators, but they might have a more
geometric description.  Since the operator in question is associated
to $B$ field holonomies on exceptional divisors,
sometimes people describe this phenomenon as having a 
``$B$ field at a quotient singularity'' -- if taken literally,
what would such a statement mean?
Can this be understood mathematically, or is it an intrinsically-CFT
phenomenon?  Why, directly in the language of the CFT
(as opposed to a massive theory), does this have the effect
of making the CFT well-behaved?  And perhaps best of all,
how can one see such a $B$ field directly in the string orbifold
CFT?  (The arguments given in \cite{edstrings95,paulz2,kw} and elsewhere
do not directly speak to the CFT itself, but either to deformations
of it, or massive theories believed to flow to it in the IR.)

\item ``String orbifold CFT's are well-behaved because they describe
strings on (`infinitesimal') resolutions.'' 

Since string orbifold CFT's are well-behaved ({\it i.e.}, they behave
as though they described strings on smooth spaces), and since they
have twist fields which are often associated to exceptional divisors
in a minimal resolution of the underlying quotient space,
physicists sometimes speak loosely of string orbifold CFT's as describing
strings on resolutions.
Of course, this cannot be literally correct in general,
simply because terminal singularities exist
 -- not all quotient singularities of interest to physics can be
resolved (and even those that can, typically do not have unique resolutions).  
Thus, any attempt to describe properties of string
orbifolds in terms of a resolution simply cannot possibly be
the general story.

Having said that, given an orbifold that can be
resolved,
it is true that the full string theory
(not just the CFT) will often see that resolution -- as fields in the
full string theory fluctuate, twist fields corresponding to blowup modes
will be excited, and so the full string theory will probe resolutions
as it probes other small deformations of the original CFT.
So, again, describing string orbifold CFT's in terms of strings
on resolutions gives a not unreasonable intuition for some of
their properties.
However, the fact that the full string theory is well-behaved
close to a given point in string moduli space is hardly itself
evidence that the CFT at that point is well-behaved -- if it were,
badly-behaved CFT's would be far more rare, as one can usually deform
them to better-behaved CFT's.
\end{itemize}

Stacks, on the other hand, seem to give a different perspective,
which appears to be simpler and cleaner.
Quotient stacks, the target spaces
of string orbifolds viewed as sigma models, are smooth\footnote{
Technically, $[X/G]$ is smooth if $X$ is smooth and $G$ is a 
discrete group acting on $X$ by diffeomorphisms.},
and, moreover, smooth in precisely the sense relevant for sigma models. 

The fact that quotient stacks themselves are smooth is not itself sufficient
to explain why string orbifold CFT's are smooth, one must also check
whether the sense in which they are smooth is physically relevant.
In this particular case, the sense in which a quotient stack is
smooth precisely coincides with the notion of smoothness relevant
for a sigma model to be well-behaved.  

However, the bottom line is that one can now see that string orbifold
CFT's are well-behaved because they are sigma models on smooth spaces,
namely, quotient stacks.

The business involving the $B$ field also appears to have
a new understanding from a stack
perspective.
To see how it arises, we shall consider D-brane
probes, described as coherent sheaves.  For readers who might be
leery of this usage in a stack context, we should point out two 
important facts:
\begin{itemize}
\item First, a coherent sheaf on a quotient stack $[X/G]$ is precisely
a $G$-equivariant coherent sheaf on $X$ (which is not quite the
same as a sheaf on the quotient space $X/G$).  Recall that the
Douglas-Moore construction \cite{dougmoore} of D-branes
on string orbifolds describes $G$-equivariant
objects on the covering space, so in other words, the Douglas-Moore
construction of branes on orbifolds precisely corresponds to
coherent sheaves on quotient stacks.
\item Second, a coherent sheaf on a flat gerbe is the same
thing as a `twisted' sheaf on the underlying space, {\it i.e.},
twisted in the sense of `bundles' on D-branes with $B$ fields.
Put another way, sheaves on gerbes are an alternative to 
modules over Azumaya algebras as popularized in \cite{kapustin}.
\end{itemize}
Given that all stacks look locally like either orbifolds or gerbes,
these two cases justify using coherent sheaves to describe
D-branes on more general stacks.

Now, in order to see what quotients stacks have to do with $B$ fields,
let us consider a naive `blowup' of the stack $[{\bf C}^2 / {\bf Z}_2 ]$.
In particular, the minimal resolution of the quotient
singularity ${\bf C}^2/{\bf Z}_2$ is the same as the quotient
$( \mbox{Bl}_1 {\bf C}^2 ) / {\bf Z}_2$, where the ${\bf Z}_2$
action has been extended trivially over the exceptional divisor
of the blowup.  Thus, the quotient stack 
$[ ( \mbox{Bl}_1 {\bf C}^2 ) / {\bf Z}_2 ]$
is a naive stacky analogue of the resolution of the quotient space
${\bf C}^2/{\bf Z}_2$, and among other things, is a stack over
the resolution.

Finally, consider D-brane probes of this stack, viewed
as coherent sheaves.
Away from the exceptional divisor, this stack looks like the 
corresponding space, so a D-brane away from the exceptional divisor
thinks it is propagating on the underlying space.
A coherent sheaf over the exceptional divisor, on the other hand,
is a sheaf on a gerbe, and so describes a D-brane in the presence
of a $B$ field.  Thus, we see the advertised $B$ field.
In fact, we can read off even more -- the gerbe over the exceptional
divisor is a ${\bf Z}_2$-gerbe, so the corresponding $B$ field holonomy
must be either $0$ or $1/2$.  Determining which requires a detailed
calculation, but notice that
without doing any work, we have quickly arrived in the right
ballpark.

Although we have spoken about D-brane probes, the same remarks
also hold true for sigma models, using the result that
a sigma model on a flat gerbe is equivalent to a sigma model
on the underlying space with a nontrivial $B$ field.

Having said all of this, in order to properly understand
the old lore concerning $B$ fields at quotient singularities would require
a detailed understanding of K\"ahler moduli of quotient stacks,
something that does not presently seem to exist.
However, the point is that from the perspective of quotient stacks,
the lore concerning $B$ fields is extremely natural -- having
a ``$B$ field at a quotient singularity'' is a natural notion forced
on one by virtue of working with stacks, whereas understanding
such a statement in terms of quotient spaces seems impossibly obscure.

A few more words should be said.  We have argued that
quotient stacks contain within themselves an intrinsic notion of
a ``$B$ field at a singularity,'' and thereby potentially cleared
up one confusing issue involving string orbifolds.  However, linking
this perspective up with traditional computations is a nontrivial matter.
Historically, physicists (implicitly assuming that the deformation
theory relevant to string orbifolds was that of quotient {\it spaces})
arrived at the same conclusions about the occurrence of a $B$ field
by considering rational curve counting arguments in limits in which
exceptional divisors shrink to zero size.  To be consistent, once we
accept that string orbifolds describe strings on quotient stacks,
any K\"ahler deformation must be a K\"ahler deformation of the {\it stack},
not the quotient space, and the resulting stack may well have
a different rational curve count.  In order for the statements that
the physics community have assumed to hold true, we need some
theorems regarding 
the relationship between rational curves in resolutions
of quotient spaces and K\"ahler deformations of quotient stacks,
something we will speak about at greater length in the next subsection.

\subsection{Deformation theory and other unresolved issues}

In the last several subsections we have described a number of
positive developments in understanding string orbifolds as
strings on quotient stacks:
\begin{itemize}
\item We have described the notion of a sigma model on a general
stack with an atlas, and verified that, indeed, within the
context of that definition a string orbifold is literally a sigma
model on a quotient stack -- in other words, at least formally
at a classical level, a string orbifold appears to be the 
same as a string compactified on a quotient stack.
\item We described the massless spectrum of a sigma model
on a quotient stack, and in so doing, resolved the apparent contradiction
that the massless spectrum of a string orbifold is {\it not} the cohomology
of the target of the purported corresponding sigma model. 
\item We have used this description of string orbifolds to
shed light on some physical properties of string orbifolds that have
historically been considered very mysterious.
\end{itemize}

However, it must be emphasized that there is still a tremendous amount
of work that must be done to completely nail down these ideas.
Even to completely nail down the notion of a sigma model on a stack,
much more work must be done to answer questions ranging from
\begin{quote}
How does one describe sigma model anomalies in the context
of stack targets?
\end{quote}
to
\begin{quote}
Sigma models with stack targets are one thing, but can one do quantum field 
theory directly on a stack,
where the stack is now taken to be the underlying spacetime?
Answering this question positively would imply being able to
explicitly write down, for example, the propagator for a scalar field
on the stack $[ {\bf R}^4 / {\bf Z}_2 ]$.
More generally, for sigma models, if the target looks like a space merely
locally, one can go a long way, but in order to do quantum field theory
on a space, global properties are important from the outset.
\end{quote}

One of the larger remaining unresolved issues involves deformation theory.
For the last fifteen years, physicists have assumed that the
deformation theory of a string orbifold is the same as that of
a quotient space.  This assumption has figured into topics ranging from
mirror symmetry to string/string duality. 
Indeed, the fact that this assumption has given consistent results
is a strong argument to many physicists
that string orbifolds describe strings on quotient {\it spaces},
that stacks have no direct physical relevance.

Clearly, one of the more important questions that must be
answered before physicists will agree that quotient stacks are physically
relevant,
is simply, why have we been able to get away with assuming the relevant
deformation theory is that of a quotient {\it space} ?
For example, K\"ahler deformations of string orbifolds have
historically been described by physicists in terms of
resolutions of the underlying quotient space $X/G$,
because physicists have made naive assumptions that string
orbifolds describe strings on quotient spaces (suitably decorated).
Once one accepts that string orbifolds describe strings on
something other than quotient spaces, these old arguments must 
be completely reexamined.  For example, a K\"ahler deformation of a 
string orbifold must be a stacky ``resolution'' of
the quotient stack $[X/G]$, if string orbifolds really do describe
strings on quotient stacks, {\it not} a resolution of the
quotient space $X/G$.  One cannot consistently speak of
string orbifolds in terms of stacks, and also simultaneously
talk about deformations and resolutions of quotient spaces.

A complete understanding of the deformation theory of quotient stacks,
even if only for a family of simple examples such as $[ {\bf C}^2 / G ]$,
would be extremely desirable.
Also, a complete understanding would probably also shed light on
other matters -- for example,
the author would not be at all surprised if yet another
way of thinking about the McKay correspondence emerged from such
considerations.

There are other matters related to deformation theory that one could
also ask.  For example, in string theory, K\"ahler moduli are paired
with ``theta angles.'' It is tempting to speculate that these theta
angles might have some sort of understanding in terms of extra data
needed to completely specify stack moduli.  

Another matter concerns the lore that ``string orbifolds have
$B$ fields at quotient singularities.''  Previously we described how one
could clearly see, at least in naive calculations, 
that $B$ fields arise very naturally in 
stack contexts.  However, to actually calculate holonomies on exceptional
divisors requires precisely
understanding the stack-y analogue of blowups, for example.
Again, knowing the deformation theory is very important.

In fact, understanding the $B$-field-at-singularities business properly
implies a nontrivial consistency check involving rational curve counting
arguments, as we implied in the last subsection.
As we described there, one of the attractive features of quotient stacks,
from a physics perspective, is that they appear to give an implicit
understanding of what it means to have a ``$B$ field at a singularity,''
as is spoken of in the physics literature.  This feature was derived
in the past by considering rational curve sums in limits in which
exceptional divisors shrunk to zero size (and so could be understood
as a limit, but was very mysterious from the direct perspective
of the string orbifold CFT).  Now, to connect the implicit understanding
outlined in the previous section with the standard derivation,
implies a nontrivial relationship.  After all, as we have stressed
here, one of the implications of the idea that string orbifolds
describe strings on quotient stacks, is that the relevant deformation
theory is that of a stack.  Physicists have long equated such deformation
theory with that of a quotient space because of some slightly naive
assumptions concerning string orbifolds, namely that they describe
strings on quotient {\it spaces} (decorated with some quantum effects
at singularities, which are typically glossed over in these discussions).
In order to agree with the usual calculations physicists perform,
rational curve counting in the stacky ``resolution'' of a stack
$[X/G]$, combined with any nontrivial $B$ field holonomy
(from any gerbe structure on exceptional divisors),
must combine to yield physical results equivalent to those
obtained by just counting rational curves in a resolution
of the quotient space $X/G$.  In other words, on the face of it,
in order to be consistent, stacky effects combined with
stacky rational curves had better combine to agree with
rational curve counting in resolutions of quotient spaces.
If this does not happen, then either string orbifolds
do not describe strings on quotient stacks, because of some
nonobvious subtlety, or the physics analysis of string moduli spaces
must be modified, which would have serious ramifications for our
understanding of everything from mirror symmetry to string/string
duality.

Yet another matter concerns intermediate points in the moduli space.
In string theory, one can interpolate between the string orbifold
point in moduli space (with `nonzero $B$ field') and the singular
quotient space point (with `zero $B$ field').
One could ask, how can the intermediate points be interpreted?
Is there a family of stacks interpolating between quotient stacks
and quotient spaces, {\it i.e.}, do all of those intermediate points
have stack interpretations?  If not, perhaps some slight generalization
of stacks is required to understand those points, or perhaps those
points can only be understood within conformal field theory.
The usual physics picture is that walking along those intermediate
points corresponds to changing the vacuum expectation value of
some operator describing strings on a (decorated) quotient space.
To many physicists, that picture sounds much more natural than the
idea that an operator which previously had no geometric interpretation,
acts to interpolate between distinct `spaces,' {\it i.e.},
between quotient spaces and quotient stacks.  This also has been
used as an indirect argument in some quarters that 
string orbifolds are not the same as strings on quotient stacks.

\section{Conclusions}

In this lecture we have given an overview of two recent developments
tied to string orbifolds in physics.

First, we gave a mostly complete understanding of discrete torsion,
a degree of freedom associated with string orbifolds, simply as the
choice of orbifold group action on a field in the theory known as the
$B$ field.  We outlined how to derive the $H^2(G, U(1))$ classification,
Vafa's twisted sector phase factors, Douglas's projectivized equivariant
K-theory, and analogues for other tensor field potentials.
Some work remains to be done to completely understand supersymmetric
moduli obstruction.

Second, we described the beginnings of a program to understand
the role of quotient stacks in string orbifolds.
Although someone not acquainted with the physics lore
might claim they `know' that
string orbifolds are the same as strings compactified on quotient stacks,
there is a tremendous amount of work that must be done to
be able to justify such a statement -- not only to make sense
out of the notion of a string compactified on {\it any} stack,
but also to reconcile such claims with the role that string orbifolds
play in the rest of string theory.  We have outlined the basics, such as
how one defines a sigma model, and how to resolve some of the basic
apparent paradoxes, such as the fact that the massless spectrum of
such sigma models is not given by the cohomology of the target space,
but a tremendous amount of work remains to be done,
In particular, this issue appears to have
nontrivial implications:
\begin{itemize}
\item If the physics lore is correct, and string orbifolds do not
describe strings on quotient stacks (and indeed, we have described
several arguments that support such a conclusion),
then some physicists would question whether a description of 
twist fields in terms of group actions on covers can hope to be
the final word on the matter.
\item If, on the other hand, string orbifolds do describe
strings on quotient spaces, then some work must be done to
properly understand string moduli spaces.  If string orbifolds
describe strings on quotient stacks, then one cannot identify
moduli of the CFT with deformations or resolutions of the
corresponding quotient {\it space}, but rather moduli of the CFT
must be identified with deformations and analogues of resolutions
for the quotient stack.
\end{itemize}
As a description of string orbifolds in terms of sigma models
on quotient stacks appears to greatly clarify the
physics of string orbifolds, we think it very important that these
issues be resolved.

\section{Acknowledgements}

The work outlined herein has been done gradually over
a long period, and so we have spoken about it to a large
number of people.  In particular, we would like to thank
P.~Aspinwall, D.~Benzvi, S.~Dean,
E.~Diaconescu, J.~Distler, M.~Douglas, T.~Gomez, J.~Harvey, S.~Hellerman, 
S.~Kachru, S.~Katz, A.~Knutson, 
G.~Moore, D.~Morrison, A.~Lawrence, T.~Pantev, R.~Plesser, 
B.~Pardon, Y.~Ruan, S.~Sethi, E.~Silverstein, M.~Stern, W.~Wang,
K.~Wendland,
and E.~Witten
for useful conversations, and we apologize to anyone whose name
we have accidentally omitted.

\end{document}